\documentclass[12pt,a4paper]{article}
\usepackage{amsmath,amsthm}
\usepackage{authblk}
\usepackage{aligned-overset}
\usepackage{newtxtext,newtxmath}
\usepackage[style=numeric-comp,maxbibnames=99,bibencoding=utf8,giveninits=true]{biblatex}
\usepackage{cancel}
\usepackage{hyperref}
\usepackage{newtxtext}
\usepackage{newtxmath}
\usepackage{graphicx}

\bibliography{handouts}

\newcommand{\email}[1]{\href{mailto:#1}{#1}}

\theoremstyle{plain}
\newtheorem{theorem}{Theorem}
\newtheorem{proposition}[theorem]{Proposition}
\newtheorem{lemma}[theorem]{Lemma}
\newtheorem{corollary}[theorem]{Corollary}

\theoremstyle{remark}
\newtheorem{remark}[theorem]{Remark}

\theoremstyle{definition}
\newtheorem{assumption}{Assumption}

\newcommand{\Real}{\mathbb{R}}

\newcommand{\Th}{\mathcal{T}_h}
\newcommand{\Fh}{\mathcal{F}_h}
\newcommand{\Fhb}{\mathcal{F}_h^{\rm b}}
\newcommand{\FT}{\mathcal{F}_T}

\newcommand{\Hdiv}[1]{H(\operatorname{div};#1)}

\newcommand{\Poly}{\mathcal{P}}
\newcommand{\RTN}{\mathcal{RT\!N}}

\newcommand{\lproj}[2]{\pi_{#1}^{#2}}
\newcommand{\eproj}[2]{\varpi_{#1}^{#2}}

\newcommand{\norm}[2]{\|#1\|_{#2}}
\newcommand{\seminorm}[2]{|#1|_{#2}}

\newcommand{\ul}[1]{\underline{#1}}

\graphicspath{{figures/}}

\newcommand{\term}{\mathfrak{T}}

\begin{document}

\title{From Finite Elements to Hybrid High-Order methods}
\author[1]{Daniele A. Di Pietro}
\author[1,2]{J\'{e}r\^{o}me Droniou}
\affil[1]{%
  IMAG, Univ Montpellier, CNRS, Montpellier 34090, France \\
  \email{daniele.di-pietro@umontpellier.fr},
  \email{jerome.droniou@umontpellier.fr}
}
\affil[2]{%
  School of Mathematics, Monash University, Melbourne, Australia
}

\maketitle

\begin{abstract}
This document contains lecture notes from the Ph.D. course given at Scuola Superiore Meridionale by Daniele Di Pietro in February 2025.
The goal of the course is to provide an overview of polytopal methods, focusing on the Hybrid High-Order (HHO) method.
As a starting point, we study the Crouzeix--Raviart method for a pure diffusion equation, with particular focus on its stability.
We then show that, switching to a fully discrete point of view, it is possible to generalize it first to polyhedral meshes and then to arbitrary order, leading to a method that belongs to the HHO family.
A study of the stability and consistency of this method reveals the need for a stabilization term, for which we identify two key properties.
\end{abstract}

\section{The Poisson problem}

Denote by $\Omega \subset \mathbb{R}^d$, $d \in \{2, 3\}$, an open polytopal (i.e., polygonal if $d = 2$ or polyhedral if $d = 3$) domain with boundary $\partial \Omega$  (the following discussion also holds for higher space dimensions, but this is not very important at this stage).
To fix the ideas, we will focus on the Poisson problem with homogeneous Dirichlet boundary conditions:
Given $f : \Omega \to \mathbb{R}$, find $u : \Omega \to \mathbb{R}$  such that
\begin{equation}\label{eq:poisson:strong}
  \begin{alignedat}{4}
    - \Delta u &= f &\qquad& \text{in $\Omega$},
    \\
    u &= 0 &\qquad& \text{on $\partial \Omega$}.
  \end{alignedat}
\end{equation}
Assuming $f \in L^2(\Omega)$, a classical weak formulation of this problem reads:
Find $u \in H_0^1(\Omega)$ such that
\begin{equation}\label{eq:poisson:weak}
  a(u, v) \coloneq \int_\Omega \nabla u \cdot \nabla v
  = \int_\Omega f v
  \qquad\forall v \in H_0^1(\Omega).
\end{equation}
The well-posedness of this weak formulation hinges on the Poincar\'e inequality stated in the following theorem.
For a proof, see, e.g., \cite[Chapter~9]{Brezis:11}.

\begin{theorem}[Continuous Poincar\'e inequality]\label{thm:poincaré:continu}
  There exists a real number $C_\Omega > 0$ only depending on $\Omega$ such that, for all $v \in H_0^1(\Omega)$,
  \begin{equation}\label{eq:poincare}
    \norm{v}{L^2(\Omega)}
    \le C_\Omega \norm{\nabla v}{L^2(\Omega)^d}.
  \end{equation}
\end{theorem}


\section{On the stability of Galerkin approximations}

Before tackling the question of the non-conforming Crouzeix--Raviart approximation of \eqref{eq:poisson:weak}, we briefly discuss
the stability of conforming Galerkin approximations for this model.

A Galerkin approximation is obtained by selecting a finite-dimensional subspace $W_{h,0}$ of $H^1_0(\Omega)$ and by writing the weak
formulation in that particular subspace: Find $u_h\in W_{h,0}$ such that
\begin{equation}\label{eq:poisson:galerkin}
  a(u_h, v_h) = \int_\Omega f v_h
  \qquad\forall v \in W_{h,0}.
\end{equation}

The \emph{stability} of a numerical approximation expresses its ability to control the amplification of errors on the data (in suitable norms): a small perturbation on the data should result in a small perturbation of the solution. For a linear model, this reduces to an a priori bound of the solution in terms of the data which, for \eqref{eq:poisson:galerkin}, is contained in the following proposition.
Its simple proof is given to illustrate the importance of the (continuous) Poincar\'e inequality in this context.

\begin{proposition}[Stability of the Galerkin approximation]\label{prop:galerkin:stability}
If $u_h\in W_{h,0}$ is a solution of \eqref{eq:poisson:galerkin}, then 
\[
\norm{\nabla u_h}{L^2(\Omega)^d}\le C_\Omega \norm{f}{L^2(\Omega)},
\]
where $C_\Omega$ is the Poincar\'e constant of Theorem \ref{thm:poincaré:continu}.
\end{proposition}

\begin{proof}
We choose $v_h=u_h$ as a test function in \eqref{eq:poisson:galerkin} and use the Cauchy--Schwarz inequality on the right-hand side to get
\[
\norm{\nabla u_h}{L^2(\Omega)^d}^2=\int_\Omega f u_h\le \norm{f}{L^2(\Omega)}\norm{u_h}{L^2(\Omega)}\le
C_\Omega \norm{f}{L^2(\Omega)}\norm{\nabla u_h}{L^2(\Omega)^d},
\]
where the conclusion follows by applying the Poincar\'e inequality \eqref{eq:poincare} to $v=u_h$, which is a valid choice since $u_h\in W_{h,0}\subset H^1_0(\Omega)$.
Simplifying the above estimate by $\norm{\nabla u_h}{L^2(\Omega)^d}$ concludes the proof.
\end{proof}


\section{A non-conforming scheme for the Poisson problem on conforming simplicial meshes}\label{sec:cr}

\subsection{The Crouzeix--Raviart space}

Denote by $T \subset \mathbb{R}^d$ a $d$-simplex, i.e., a triangle if $d=2$ or a tetrahedron if $d=3$, and by $\FT$ the set collecting its faces.
For the sake of simplicity, from this point on we will only use the three-dimensional terminology, i.e., we will speak of polyhedra and faces instead of polygons and edges.

Adopting the classical definition of \cite{Ciarlet:02}, the Crouzeix--Raviart element, originally introduced in \cite{Crouzeix.Raviart:73} for the discretisation of incompressible flows\footnote{In \cite{Crouzeix.Raviart:73}, the authors actually introduce a velocity--pressure element. For this reason, some authors refer to the scalar  element considered here as ``non-conforming $\mathbb{P}^1$'' rather than ``Crouzeix--Raviart''.}, is the triplet $(T, \Poly^1(T), \sigma_T)$, where $\Poly^1(T)$ is spanned by the restriction to $T$ of affine functions, while the degrees of freedom $\sigma_T = ( \sigma_{TF} )_{F \in \FT}$ are such that, for all $F \in \FT$,
\[
\sigma_{TF} : \Poly^1(T) \ni v \mapsto v_F \coloneq \frac{1}{|F|} \int_F v \in \mathbb{R}.
\]

Denote by $\Th$ a tetrahedral mesh and by $\Fh$ the set collecting its faces.
We will additionally assume, as customary in finite elements, that $\Th$ is \emph{conforming}, i.e., that the intersection of two distinct mesh elements $T_1$ and $T_2$ of $\Th$ is either empty or coincides with a mesh vertex, edge, or, if $d = 3$, a face.
We consider the global Crouzeix--Raviart space $V_h$ spanned by piecewise affine functions $v_h$ that have continuous average value across interfaces, i.e., such that $\sigma_{T_1F} (v_{h|T_1}) = \sigma_{T_2F} (v_{h|T_2})$ for every $F \in \Fh \setminus \Fhb$ shared by the mesh elements $T_1$ and $T_2$; here, $\Fhb \coloneq \left\{ F \in \Fh \,:\, F \subset \partial \Omega \right\}$ is the set of boundary faces.
To account for boundary conditions, we consider the following subspace of $V_h$:
\[
V_{h,0} \coloneq \left\{
v_h \in V_h \,:\, \text{$v_F = 0$ for all $F \in \Fhb$}
\right\}.
\]

\subsection{Broken gradient}

The space $V_h$ is a subspace of the broken Sobolev space
\[
H^1(\Th) \coloneq \left\{ v \in L^2(\Omega) \,:\, \text{$v_{|T} \in H^1(T)$ for all $T \in \Th$} \right\}
\]
but not of $H^1(\Omega)$.
Therefore, we cannot apply the standard gradient to the elements of $V_h$.
We can, however, apply the broken gradient $\nabla_h : H^1(\Th) \to L^2(\Omega)^d$ such that, for all $v \in H^1(\Th)$, $(\nabla_h v_h)_{|T} \coloneq \nabla v_{h|T}$.

\subsection{A non-conforming scheme for the Poisson problem}\label{sec:cr:discrete.problem}

Replacing the Sobolev space $H^1_0(\Omega)$ with the discrete space $V_{h,0}$ and the standard gradient with the broken gradient leads to the following discrete problem:
Find $u_h \in V_{h,0}$ such that
\begin{equation}\label{eq:cr:discrete}
a_h(u_h, v_h) \coloneq \int_\Omega \nabla_h u_h \cdot \nabla_h v_h
= \int_\Omega f v_h\qquad\forall v_h\in V_{h,0}.
\end{equation}

We will assume, from this point on, that $\Th$ belongs to a refined sequence of finite element meshes that is regular in the usual sense of \cite{Ciarlet:02}.
We will use the symbol $\lesssim$ for inequalities that hold up to a positive constant independent of the meshsize and, for local inequalities on a mesh element or face $Y$, also on $Y$, but possibly depending on other quantities such as the domain $\Omega$, the dimension $d$, and the mesh regularity parameter.

Since $V_{h,0} \not \subset H_0^1(\Omega)$, the continuous Poincar\'e inequality \eqref{eq:poincare} is not inherited at the discrete level.
As a consequence, the stability of the scheme cannot be directly proved as for conforming Galerkin approximations (see the proof of Proposition \ref{prop:galerkin:stability}). Instead, we need to develop the following discrete counterpart of this inequality, valid in the Crouzeix--Raviart space.

\begin{lemma}[Discrete Poincar\'e inequality for the Crouzeix--Raviart space]\label{lem:poincare:cr}
  For all $v_h \in V_{h,0}$, it holds
  \[
  \norm{v_h}{L^2(\Omega)} \lesssim \norm{\nabla_h v_h}{L^2(\Omega)^d}.
  \]
\end{lemma}

The proof of this lemma will be given in Section \ref{sec:proof:poincare:cr} after recalling some relevant results.
The stability of the Crouzeix--Raviart scheme then follows exactly as in the proof of Proposition \ref{prop:galerkin:stability}, substituting the discrete Poincar\'e inequality above for the continuous one. We therefore omit the proof of the following corollary.

\begin{corollary}[Stability of the Crouzeix--Raviart scheme]
If $u_h\in V_{h,0}$ is a solution of \eqref{eq:cr:discrete}, then
\[
\norm{\nabla_hu_h}{L^2(\Omega)^d}\lesssim \norm{f}{L^2(\Omega)}.
\]
\end{corollary}

\subsection{A magic formula}

A relevant space in the study of the Poisson problem is
\[
\Hdiv{\Omega} \coloneq \left\{
\tau \in L^2(\Omega)^d \,:\, \operatorname{div} \tau \in L^2(\Omega)
\right\}
\]
equipped with the norm
\[
\norm{\tau}{H(\operatorname{div}; \Omega)}
\coloneq \left(
\norm{\tau}{L^2(\Omega)^d}^2
+ \norm{\operatorname{div} \tau}{L^2(\Omega)}^2
\right)^{\frac12}.
\]
A crucial property of this space is that its elements have continuous normal trace across interfaces.
Specifically, let $\tau \in \Hdiv{\Omega}$ and further assume that $\tau$ is sufficiently regular inside each element $T \in \Th$ for its normal trace to be in $L^2(F)$ on each face $F \in \FT$ (this is the case, e.g., if $\tau \in H^1(\Th)^d$).
Then, for each mesh interface $F \in \Fh \setminus \Fhb$ shared by distinct elements $T_1$ and $T_2$, it holds
\begin{equation}\label{eq:continuity.normal.components}
  \tau_{|T_1} \cdot n_{T_1F} + \tau_{|T_2} \cdot n_{T_2F} = 0\quad\mbox{on $F$},
\end{equation}
where, for any $T\in\Th$ and $F \in \FT$, $n_{TF}$ denotes the unit vector normal to $F$ and pointing out of $T$.
This property is at the heart of the following magic formula.

\begin{proposition}[Magic formula]\label{prop:magic.formula}
  Denote by $( v_F )_{F \in \Fh}$ a family of functions $v_F \in L^2(F)$ such that $v_F = 0$ for all $F \in \Fhb$ and let $\tau \in H(\operatorname{div}; \Omega)$ be such that, for all $T\in\Th$ and $F\in\mathcal F_T$, the normal trace of $\tau_{|T}$ on $F$ is in $L^2(F)$.
  Then,
  \begin{equation}\label{eq:magic.formula}
    \sum_{T \in \Th} \sum_{F \in \FT} \int_F (\tau_{|T} \cdot n_{TF}) v_F = 0.
  \end{equation}
\end{proposition}

\begin{proof}
  For all $F \in \Fh$, denote by $\mathcal{T}_F$ the set containing the one (if $F$ is a boundary face) or two (if $F$ is an interface) elements that have $F$ as a face.
  Exchanging the sums in the left-hand side of \eqref{eq:magic.formula}, we obtain
  \begin{align}
    &\sum_{T \in \Th} \sum_{F \in \FT} \int_F (\tau_{|T} \cdot n_{TF}) v_F
    \nonumber\\
    &\quad
    = \sum_{F \in \Fh} \sum_{T \in \mathcal{T}_F} \int_F (\tau_{|T} \cdot n_{TF}) v_F
    \nonumber\\    
    &\quad
    = \sum_{F \in \Fh \setminus \Fhb} \int_F\left(\sum_{T\in\mathcal T_F}\tau_{|T}\cdot n_{TF}\right)v_F
    + \sum_{F \in \Fhb}  \int_F (\tau_{|T_F} \cdot n_{T_F F}) v_F,
    \label{eq:magic.formula.swap}
  \end{align}
  where $\mathcal{T}_F = \{ T_F \}$ if $F \in \Fhb$.
  Noticing that \eqref{eq:continuity.normal.components} precisely states that, for all $F\in\mathcal F_h\backslash \mathcal F_h^{\rm b}$, $\sum_{T\in\mathcal T_F}\tau_{|T}\cdot n_{TF}=0$, and recalling that $v_F=0$ whenever $F \in \Fhb$, \eqref{eq:magic.formula.swap} concludes the proof of \eqref{eq:magic.formula}.
\end{proof}

\begin{remark}[Validity on general meshes]\label{rem:magic.formula:general.meshes}
  The proof of Proposition~\ref{prop:magic.formula} does not use the fact that $\Th$ is a matching simplicial mesh.
  As a matter of fact, this result holds true for general decompositions of $\Omega$ into sets with a piecewise smooth boundary.
\end{remark}

\subsection{The Raviart--Thomas--N\'ed\'elec space}

The second tool that we need to prove Lemma \ref{lem:poincare:cr} is the Raviart--Thomas--N\'ed\'elec space \cite{Raviart.Thomas:77,Nedelec:80}.
Let $T \in \Th$ and set
\begin{equation}\label{eq:RTN}
  \RTN^1(T) \coloneq \Poly^0(T)^d + (x - x_T) \Poly^0(T),
\end{equation}
where $\Poly^0(T)$ is spanned by constant functions on $T$ and $x_T$ is a point inside $T$.
The space \eqref{eq:RTN} is a \emph{trimmed} polynomial space in the sense that it sits between $\Poly^0(T)^d$ and $\Poly^1(T)^d$, without coinciding with either of these spaces.
Letting $\sigma_T = ( \sigma_{TF} )_{F \in \FT}$ be such that
\[
\sigma_{TF} : \RTN^1(T) \ni \tau \mapsto \frac{1}{|F|} \int_F \tau \cdot n_{TF} \in \mathbb{R},
\]
we obtain the finite element $(T, \RTN^1(T), \sigma_T)$.

\begin{proposition}[Divergence and normal trace of functions in $\RTN^1(T)$]\label{prop:div.normal.RTN}
	For all $T \in \Th$ and all $\tau \in \RTN^1(T)$,
  \begin{equation}\label{eq:div.trace.RTN1}
    \text{%
      $\operatorname{div} \tau \in \Poly^0(T)$
and $\tau \cdot n_{TF} \in \Poly^0(F)$ for all $F \in \FT$.      
    }
  \end{equation}
\end{proposition}

\begin{proof}  
  By definition \eqref{eq:RTN} of $\RTN^1(T)$, there exists $\overline{\tau} \in \Poly^0(T)^d$ and $q \in \Poly^0(T)$ such that $\tau = \overline{\tau} + (x - x_T) q$.
  Clearly, $\operatorname{div} \tau = d q \in \Poly^0(T)$.
  Let now $F \in \FT$, and notice that, for all $x \in F$, $(x - x_T) \cdot n_{TF}\equiv \operatorname{dist}(x_T,H_F)$ is the distance between $x_T$ and the hyperplane $H_F$ spanned by $F$.
  Hence, $\tau \cdot n_{TF} = \overline{\tau} \cdot n_{TF} + \operatorname{dist}(x_T, F) \in \Poly^0(F)$ on $F$, thus concluding the proof.
\end{proof}

The global Raviart--Thomas--N\'ed\'elec space is obtained setting
\[
\Sigma_h \coloneq \left\{
\tau \in H(\operatorname{div}; \Omega) \,:\, \text{$\tau_{|T} \in \RTN^1(T)$ for all $T \in \Th$}
\right\}.
\]
The following proposition is a special case of classical results for arbitrary-order Raviart--Thomas--Nédélec spaces; see, e.g., \cite[Proposition~2.3.3]{Boffi.Brezzi.ea:13} or \cite[Lemma~3.6]{Gatica:14}.

\begin{proposition}[Surjectivity of the divergence from $\RTN^1(\Th)$ onto $\Poly^0(\Th)$]\label{prop:surjectivity.div.RTN1.P0}
  Denote by $\Poly^0(\Th)$ the space of piecewise constant functions on the mesh.
  Then, for all $v_h \in \Poly^0(\Th)$, there exists $\tau_h \in \RTN^1(\Th)$ such that
  \begin{equation}\label{eq:surjectivity.div.RTN1.P0}
    \text{%
      $\operatorname{div} \tau_h = v_h$ and
      $\norm{\tau_h}{\Hdiv{\Omega}} \lesssim \norm{v_h}{L^2(\Omega)}$.
    }
  \end{equation}
\end{proposition}

\begin{proof}
  Since $v_h \in L^2(\Omega)$, by the surjectivity of $\operatorname{div} : H^1(\Omega)^d \to L^2(\Omega)$ (see, e.g. \cite{Girault.Raviart:86,Bogovskii:80,Solonnikov:01,Duran.Muschietti:01}), there exists $\tau \in H^1(\Omega)^d$ such that
  \begin{equation}\label{eq:surjectivity.div.H1.L2}
    \text{
      $\operatorname{div} \tau = v_h$
      and $\norm{\tau}{H^1(\Omega)^d} \lesssim \norm{v_h}{L^2(\Omega)}$.
    }
  \end{equation}
  Denote by $\tau_h$ the interpolate of $\tau$ on $\Sigma_h$ such that
  \begin{equation}\label{eq:I.RTN}
    \forall T \in \Th,\qquad
    \int_F \tau_h \cdot n_{TF} = \int_F \tau \cdot n_{TF}
    \qquad \forall F \in \FT.
  \end{equation}
  By the boundedness of the Raviart--Thomas--Nédélec interpolator proved, e.g., in \cite[Lemma~4.4]{Gatica:14}, it holds
  \[
  \norm{\tau_h}{\Hdiv{\Omega}}
  \lesssim \norm{\tau}{H^1(\Omega)^d}
  \overset{\eqref{eq:surjectivity.div.H1.L2}}\lesssim \norm{v_h}{L^2(\Omega)},
  \]
  which is the second condition in \eqref{eq:surjectivity.div.RTN1.P0}.
  Moreover, for all $T \in \Th$,
  \[
  \begin{aligned}
    (v_h)_{|T}
    \equiv \frac{1}{|T|} \int_T \operatorname{div} \tau
    &= \frac{1}{|T|} \sum_{F \in \FT} \int_F \tau \cdot n_{TF}
    \\
    \overset{\eqref{eq:I.RTN}}&=  \frac{1}{|T|} \sum_{F \in \FT} \int_F \tau_h \cdot n_{TF}
    = \frac{1}{|T|} \int_T \operatorname{div} \tau_h
    \overset{\eqref{eq:div.trace.RTN1}}\equiv \operatorname{div} \tau_h,
  \end{aligned}
  \]
  where we have identified constant functions in $T$ with their value whenever needed and the second and fourth equalities result from integrations by parts.
  This shows that
  \[
  v_h = \operatorname{div} \tau_h,
  \]
  which is precisely the first condition in \eqref{eq:surjectivity.div.RTN1.P0}.
\end{proof}

\subsection{Proof of Lemma \ref{lem:poincare:cr}}\label{sec:proof:poincare:cr}

We are now ready to prove the discrete Poincar\'e inequality in the Crouzeix--Raviart space.

\begin{proof}[Lemma \ref{lem:poincare:cr}]
  Let $v_h \in V_{h,0}$, set
  \[
  v_F \coloneq \frac{1}{|F|} \int_F v_h \qquad \forall F \in \Fh,
  \]
	and let $\overline{v}_h \in L^2(\Omega)$ be the piecewise constant function on $\Th$ such that
  \begin{equation}\label{eq:poincare:cr:overline.vh}
    \overline{v}_T
    \coloneq (\overline{v}_h)_{|T} 
    = \frac{1}{|T|} \int_T v_h \qquad T \in \Th.
  \end{equation}
  By Proposition~\ref{prop:surjectivity.div.RTN1.P0}, there exists $\tau_h \in \Sigma_h$ such that
  \begin{equation}\label{eq:surjectivity.div.RTN.P0}
    \text{%
      $\operatorname{div} \tau_h = \overline{v}_h$
      and $\norm{\tau_h}{\Hdiv{\Omega}} \lesssim \norm{\overline{v}_h}{L^2(\Omega)}$.
    }
  \end{equation}
  We can then write
	\begin{equation}\label{eq:poincare:cr:basic}
  \begin{aligned}
    \norm{\overline{v}_h}{L^2(\Omega)}^2  
    \overset{\eqref{eq:surjectivity.div.RTN.P0}}&=
    \int_\Omega  \overline{v}_h \operatorname{div} \tau_h
    \\
    &= \sum_{T \in \Th} \int_T v_h \operatorname{div} \tau_h
    \\
    &=  -\int_\Omega \nabla_h v_h \cdot \tau_h
    + \sum_{T \in \Th} \sum_{F \in \Fh} \int_F v_F (\tau_h \cdot n_{TF}),
  \end{aligned}
	\end{equation}
  where the passage to the second line is justified observing that, for all $T \in \Th$, since $(\operatorname{div} \tau_h)_{|T}$ is constant inside $T$ by~\eqref{eq:div.trace.RTN1}, and identifying constant functions with their value,
  \[
  \int_T \overline{v}_T\operatorname{div} \tau_h
  = |T| \overline{v}_T (\operatorname{div} \tau_h)_{|T}
  \overset{\eqref{eq:poincare:cr:overline.vh}}=
  \left(\int_T v_h\right) (\operatorname{div} \tau_h)_{|T}
  = \int_T v_h\operatorname{div} \tau_h.
  \]
	We next notice that $\tau_h$ and $( v_F )_{F \in \Fh}$ meet the assumptions of Proposition \ref{prop:magic.formula}, hence the second term in the right-hand side of \eqref{eq:poincare:cr:basic} is zero.
  Using a Cauchy--Schwarz inequality, we can then go on writing
  \[
  \norm{\overline{v}_h}{L^2(\Omega)}^2
  \le \norm{\nabla_h v_h}{L^2(\Omega)^d} \norm{\tau_h}{L^2(\Omega)^d}
  \overset{\eqref{eq:surjectivity.div.RTN.P0}}\lesssim
  \norm{\nabla_h v_h}{L^2(\Omega)^d} \norm{\overline{v}_h}{L^2(\Omega)}.
  \]
  Simplifying, we get
  \begin{equation}\label{eq:poincare:overline.vh}
    \norm{\overline{v}_h}{L^2(\Omega)}
    \lesssim \norm{\nabla_h v_h}{L^2(\Omega)^d}.
  \end{equation}

  It only remains to control the variations of $v_h$ inside each element with the $L^2$-norm of its broken gradient.
  To this purpose, let $T \in \Th$ and set, for the sake of brevity, $v_T \coloneq (v_h)_{|T}$.
  For all $x \in T$, it holds $v_T(x) = \overline{v}_T + (x - \overline{x}_T) \cdot \nabla v_T$ with $\overline{x}_T \coloneq \frac{1}{|T|} \int_T x$ center of mass of $T$. Hence,
  \[
  \norm{v_T - \overline{v}_T}{L^2(T)}
  = \norm{(x - \overline{x}_T) \cdot \nabla v_T}{L^2(T)}
  \le h_T \norm{\nabla v_T}{L^2(T)^d},
  \]
  where the conclusion follows using a H\"older inequality followed by $\norm{x - \overline{x}_T}{L^\infty(T)^d} \le h_T$.
	We thus have
  \[
  \norm{v_T}{L^2(T)}
  \le \norm{v_T - \overline{v}_T}{L^2(T)}^2
  + \norm{\overline{v}_T}{L^2(T)}
  \lesssim h_T \norm{\nabla v_T}{L^2(T)^d}
  + \norm{\overline{v}_T}{L^2(T)}.
  \]
  Squaring the above inequality, using the fact that $(a + b)^2 \le 2 a^2 + 2 b^2$ for all real numbers $a$ and $b$, noticing that $h_T \le \operatorname{diam}(\Omega) \lesssim 1$, and summing over $T \in \Th$, we get
  \[
  \norm{v_h}{L^2(\Omega)}^2
  \lesssim \norm{\nabla_h v_h}{L^2(\Omega)^d}^2
  + \norm{\overline{v}_h}{L^2(\Omega)}^2
  \overset{\eqref{eq:poincare:overline.vh}}\lesssim
  \norm{\nabla_h v_h}{L^2(\Omega)^d}^2.
  \]
  Taking the square root of the above inequality concludes the proof.
\end{proof}


\section{A hybrid high-order space on general meshes}

\subsection{Limitations of the finite element approach}

While finite elements are more flexible than other discretization methods for PDEs, they can display limitations in practical applications.
We will discuss two examples: limitations related to the conformity constraint on the mesh, and limitations related to the need to define an underlying space of functions.

\begin{figure}\centering
  \includegraphics[height=3.25cm]{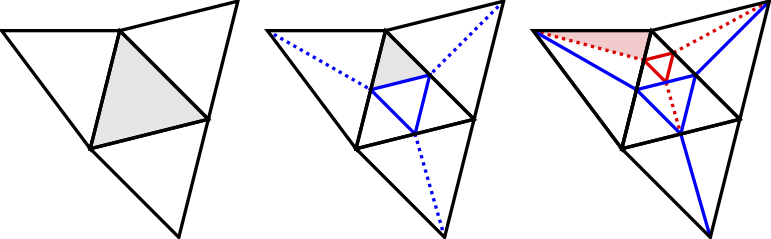}
  \caption{Two successive local refinements of a conforming finite element mesh. The element to be refined at the next step is in gray. As the refinement proceeds, and unless specific measures are taken, more and more elongated elements, such as the one in red, appear.}
  \label{fig:conforming.mesh.refinement}
\end{figure}

Finite element spaces are naturally constructed on conforming meshes.
To understand why this can be a limiting factor, let us consider the two-dimensional example in Figure~\ref{fig:conforming.mesh.refinement}.
The situation is that of an adaptive computation where, after solving the discrete problem on a given mesh and estimating the contribution to the error from each element, one refines the elements where the error is larger.
More specifically, the picture shows an initial finite element mesh and two local refinement steps.
The elements to be refined at the next step are marked in gray.
The refinement of these elements is done by joining the edge middlepoints, thus splitting the original triangle into four similar triangles with halved diameter.
In order for the mesh to remain conforming, the refinement has to be propagated to the neighbors in the sense of faces, that are subdivided into two triangles by the dotted lines.
Treating the neighbors this way prevents further propagation of the refinement.
After repeating this procedure twice, one observes the appearance of elongated elements like the one in red.
Such elements have a diameter-to-radius ratio which is much larger than the original ones, and such ratio could potentially grow larger if we repeated again and again the local refinement step.
In other words, the constraint of working with conforming meshes has forced us to \emph{trade mesh quality for resolution}.
It is, of course, possible to mitigate this phenomenon, either by identifying more efficient splittings of neighboring elements or by moving the mesh vertices after refinement.
The first strategy can limit the degradation of the mesh quality, but not prevent it altogether.
The second strategy can have unwanted side effects, e.g., when it is necessary to project the solution on the coarser mesh onto the locally refined one, since the corresponding finite element spaces are no longer in a hierarchical relation.

\begin{figure}\centering
  \includegraphics[height=3.25cm]{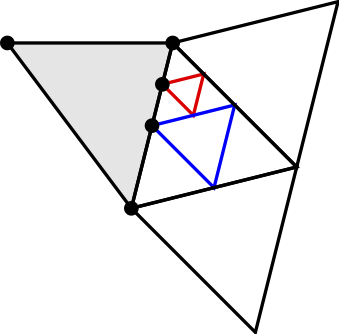}
  \caption{Non-conforming local mesh refinement. If a method supporting polygonal meshes is available, the element in gray can simply be regarded as a pentagon, and it will not require a specific treatment.}
  \label{fig:non.conforming.mesh.refinement}
\end{figure}

A more efficient way of dealing with this problem would be to remove the conformity requirement on the mesh.
Without this requirement, the local refinement process could result in meshes like the one depicted in Figure~\ref{fig:non.conforming.mesh.refinement}, where the elements marked for refinement have been subdivided as before, but the neighbors have been left untouched.
Each element of the refined mesh is similar to one element of the original mesh, so mesh quality is preserved.
While it is possible to deal with this kind of meshes in finite elements by freezing the non-conforming degrees of freedom, this requires a specific treatment and can reduce the approximation properties of the elements affected by the freezing (thus defeating the purpose of local refinement, which requires the additional degrees of freedom to be free in order to locally improve the approximation properties of the space).
If, on the other hand, we disposed of a method supporting polygonal meshes, the situation depicted in Figure~\ref{fig:non.conforming.mesh.refinement} could be dealt with seamlessly: elements such as the one in gray could simply be treated as polygons with parallel sides.

\begin{figure}\centering
  \includegraphics[height=3.25cm]{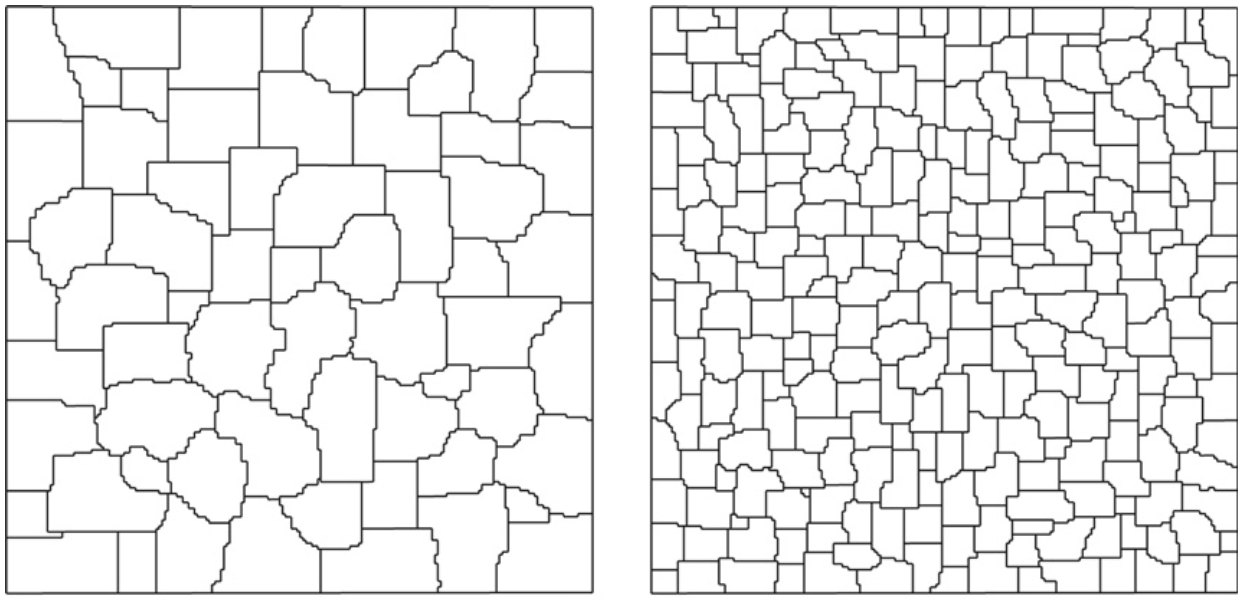}
  \caption{Two examples of meshes with, respectively, 64 and 255 elements obtained by coarsening an initial $200 \times 200$ Cartesian orthogonal mesh.}
  \label{fig:coarsened.mesh}
\end{figure}

Another limitation of finite element meshes is encountered in applications, such as the modelling of geological flows, where all the geometric information is carried by the mesh.
In such circumstances, the mesh may have elements that, from a computational point of view, are unnecessarily small in certain regions (e.g., to account for fine geometric features) or unnecessarily large in others.
While local mesh refinement is an option to solve the latter problem, it is in general not possible to coarsen a given mesh in finite element codes.
Polytopal methods, on the other hand, support much more general elements, thus paving the way to adaptive mesh coarsening.
An example of mesh coarsening is depicted in Figure~\ref{fig:coarsened.mesh}, which displays two meshes with, respectively, 64 and 255 elements obtained by coarsening an initial $200 \times 200$ Cartesian orthogonal mesh.
This idea was first highlighted in \cite{Bassi.Botti.ea:12} in the context of Discontinuous Galerkin methods; see also \cite{Antonietti.Giani.ea:13} concerning the $hp$-version of the method.

\begin{figure}
  \includegraphics[height=3.25cm]{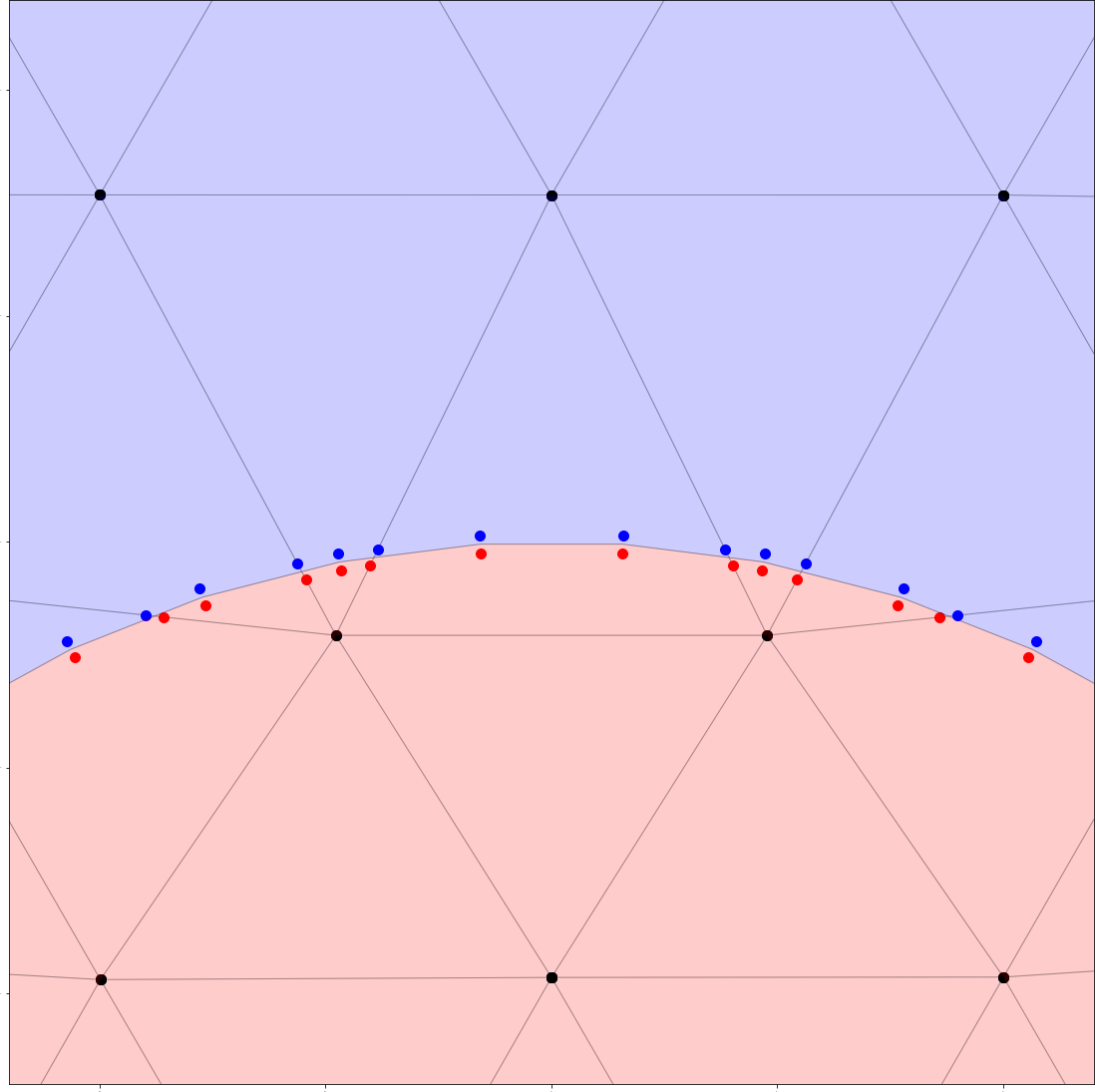}\centering
  \caption{Standard elements cut by an interface generate polygonal elements.}
  \label{fig:cut-elements}
\end{figure}

A last example of mesh-related limitation in finite element codes arises when dealing with interface problems.
In practically relevant cases involving, e.g., moving interfaces, it is not computationally feasible to generate a fitted mesh at each time step, and one has to deal with elements cut by the interface.
In the context of finite element methods, a special treatment is required for such elements.
In the Generalized/Extended Finite Element method \cite{Sukumar.Moes.ea:00,Strouboulis.Babuska.ea:00}, non-polynomial functions with compact support are added to the discrete space.
In the Immersed Finite Element method \cite{Adjerid.Babuska.ea:23}, the added functions are piecewise polynomials.
In the CutFEM method \cite{Burman.Claus.ea:15,Burman.Ern:18}, interface conditions are taken into account by using discontinuous basis functions inside the elements cut by the interface and by relying on Nitsche's techniques for their enforcement.
Polytopal methods, on the other hand, have no problem in treating each portion of an element cut by the interface as a computational element of general shape; see Figure~\ref{fig:cut-elements}.
This treatment is both seamless, in the sense that no special function or ad hoc modification of the construction is required, and robust, since, in the presence of small or elongated cuts, mesh quality can be restored by merging them into neighboring elements in the spirit of \cite{Bassi.Botti.ea:12,Antonietti.Giani.ea:13,Johansson.Larson:13}.

We conclude this section by pointing out limitations that are more subtle than the ones resulting from the restrictive notion of mesh in finite elements.
The need to identify a local function space on which the degrees of freedom are unisolvent can sometimes become a serious hurdle.
In the case of the Crouzeix--Raviart element, e.g., this constraint makes it difficult to identify higher-order versions of the space that can be used in practice; see, e.g., \cite{Ciarlet.Dunkl.ea:18}.
Another example are mixed finite elements for elasticity, where formidable challenges arise to incorporate the symmetry constraint into the space.
The first mixed finite element for elasticity with polynomial functions \cite{Arnold.Winther:02} was only discovered in 2002 , twenty-five years after mixed finite elements for vector fields were introduced \cite{Raviart.Thomas:77}.
The advent of polytopal methods was made possible by a simple but powerful idea, namely that the notion of local space could be removed altogether and replaced by reconstructions of quantities of interest directly from the degrees of freedom.
This idea was first explored in the context of lowest-order methods \cite{Kuznetsov.Lipnikov.ea:04,Bochev.Hyman:06,Codecasa.Specogna.ea:07,Beirao-da-Veiga.Lipnikov.ea:14,Droniou.Eymard:06,Droniou.Eymard.ea:10,Eymard.Gallouet.ea:10}, and then extended to the arbitrary order in a systematic way in the seminal works \cite{Beirao-da-Veiga.Brezzi.ea:13,Beirao-da-Veiga.Brezzi.ea:13*1,Di-Pietro.Ern.ea:14,Di-Pietro.Ern:15}.
As we will see in the rest of this chapter, directly working on the degrees of freedom leads to a very natural and elegant extension of the Crouzeix--Raviart methods to general meshes and arbitrary order.

\subsection{Shifting point of view}\label{sec:HHO.0:shift}

Let $v \in \Poly^1(T)$ and set
\[
v_F \coloneq \frac{1}{|F|} \int_F v\qquad
\forall F \in \FT.
\]
In what follows, we will identify $v_F$ with the constant function in $\Poly^0(F)$ taking the value $v_F$ on every point $x \in F$ whenever needed.
Since $v$ is affine, denoting for all $Y \in \{ T \} \cup \FT$ by $\overline{x}_Y \coloneq \frac{1}{|Y|} \int_Y x$ the center of mass of $Y$, it holds
\[
\frac{1}{|F|} \int_F v_F
= \frac{1}{|F|} \int_F v
= v\left( \frac{1}{|F|} \int_F x \right)
= v(\overline{x}_F).
\]
We additionally have $\nabla v \in \Poly^0(T)^d$ and, for all $w \in \Poly^1(T)$, an integration by parts gives
\begin{equation}\label{eq:cr:grad.vT}
  \int_T \nabla v \cdot \nabla w
  = -\int_T v \,\cancel{\Delta w}
  + \sum_{F \in \FT} \int_F v (\nabla w \cdot n_{TF})
  = \sum_{F \in \FT} \int_F v_F (\nabla w \cdot n_{TF}),
\end{equation}
where we have used the fact that $\Delta w = 0$ since $w$ is affine in the cancellation and the fact that $\nabla w \cdot n_{TF}$ is constant on each $F \in \FT$ (since both $\nabla w$ and $n_{TF}$ are) to replace $v$ with $v_F$ in the boundary integral.
Since $v$ is affine, $\nabla v$ is entirely determined by the left-hand side of \eqref{eq:cr:grad.vT} (for generic functions $w\in\mathcal P^1(T)$) which, in turn, can be fully computed from the averaged face values $(v_F)_{F\in\mathcal F_T}$ that appear in the right-hand side.

Let us now express the average value of $v$ in terms of the values $(v_F)_{F \in \FT}$.
We start by noticing that
\[
\begin{aligned}
  \int_T v
  &= \frac{1}{d} \int_T v \operatorname{div}(x - \overline{x}_T)
  \\
  & = - \frac{1}{d} \cancel{\int_T \nabla v \cdot (x - \overline{x}_T)}
  + \frac{1}{d} \sum_{F \in \FT} \int_F v (x - \overline{x}_T) \cdot n_{TF}
  \\
  &= \sum_{F \in \FT} \frac{d_{TF}}{d} \int_F v_F,
\end{aligned}
\]
where we have used the fact that $\operatorname{div}(x - \overline{x}_T) = d$ in the first equality,
an integration by parts in the second equality,
and noticed that $\nabla v$ is constant inside $T$ and that the function $T \ni x \mapsto x - \overline{x}_T \in \Real^d$ has zero average inside $T$ to cancel the first term in the right-hand side.
To conclude, we have noticed that
\begin{equation}\label{eq:dTF}
  (x - \overline{x}_T) \cdot n_{TF} \equiv d_{TF} \coloneq \operatorname{dist}(\overline{x}_T, H_F) \quad \forall x \in F
\end{equation}
(where, as in the proof of Proposition \ref{prop:div.normal.RTN}, $\operatorname{dist}(\overline{x}_T, H_F)$ denotes the distance of $\overline{x}_T$ from the hyperplane $H_F$ containing $F$)
and recalled the definition of $v_F$.
As a consequence
\begin{equation}\label{eq:cr:average.vT}
  \int_T v
  = \sum_{F \in \FT} \frac{d_{TF}}{d} \int_F v_F.
\end{equation}

\begin{remark}[Extension to polytopal elements]\label{rem:polytopal.extension}
  The relations \eqref{eq:cr:grad.vT} and \eqref{eq:cr:average.vT} show that a full knowledge of $v$ is provided by the face values $(v_F)_{F \in \FT}$.
  Even more interesting, \emph{the above formulas do not use the fact that $T$ is a triangle (if $d = 2$) or a tetrahedron (if $d = 3$), and actually hold for any polygon/polyhedron with planar faces that is star-shaped with respect to its center of mass}.
    The above assumptions are needed for \eqref{eq:dTF} to hold: by the planar face assumption, the quantity $(x - \overline{x}_T) \cdot n_{TF}$ is independent of $x \in F$, and the fact that $T$ is star-shaped with respect to $\overline{x}_T$ implies that this constant value is positive and equal to $d_{TF}$.
\end{remark}  

\subsection{A lowest-order hybrid space on general meshes}\label{sec:HHO:lowest.order.space}

Let $\Th$ be a polyhedral mesh and denote by $\Fh$ the corresponding set of (planar) mesh faces (as we will see in the next chapter, the notion of mesh face can differ from that of polyhedral face, but for the moment being we can safely disregard this distinction).
We define the following space
\[
\ul{V}_h^0 \coloneq \left\{
\ul{v}_h = ( v_F )_{F \in \Fh} \,:\,
\text{$v_F \in \Poly^0(F)$ for all $F \in \Fh$}
\right\}.
\]
Unlike the Crouzeix--Raviart space, which is spanned by functions $\Omega \to \mathbb{R}$, the space $\ul{V}_h^0$ is spanned by vectors of polynomials, which can be interpreted, for each face $F \in \Fh$, as the average traces of certain scalar functions on $\Omega$ which we will not bother to identify.
We shall also denote the restriction of this space to a mesh element $T \in \Th$ by $\ul{V}_T^0$, i.e.,
\[
\ul{V}_T^0 \coloneq \left\{
\ul{v}_T = ( v_F )_{F \in \FT} \,:\,
\text{$v_F \in \Poly^0(F)$ for all $F \in \FT$}
\right\}.
\]

\subsection{Affine potential reconstruction}

In light of Remark \ref{rem:polytopal.extension}, and taking inspiration from the formulas in Section \ref{sec:HHO.0:shift}, for each element $T \in \Th$ we can define an operator which, given a vector of face values in $\ul{V}_T^0$, returns an affine function mimicking the formulas \eqref{eq:cr:grad.vT} and \eqref{eq:cr:average.vT}.
Specifically, we define the \emph{potential reconstruction} $p_T^1 : \ul{V}_T^0 \to \Poly^1(T)$ such that, for all $\ul{v}_T \in \ul{V}_T^0$,
\begin{subequations}\label{eq:pT.1}
  \begin{gather}\label{eq:pT.1:gradient}
    \int_T \nabla p_T^1 \ul{v}_T \cdot \nabla w
    = \sum_{F \in \FT} \int_F v_F (\nabla w \cdot n_{TF})
    \qquad \forall w \in \Poly^1(T),
    \\ \label{eq:pT.1:average}
    \int_T p_T^1 \ul{v}_T
    = \sum_{F \in \FT} \frac{d_{TF}}{d} \int_F v_F.
  \end{gather}
\end{subequations}

\begin{remark}[Potential reconstruction]
  The naming ``potential reconstruction'' is evocative of the fact that $p_T^1$ reconstructs, from face values, a quantity that is homogeneous to $u$ in \eqref{eq:poisson:strong}, and this is, in turn, a potential in the sense of Fick's law of diffusion.
\end{remark}

\subsection{Interpolator and polynomial consistency of the affine potential reconstruction}

Let now $\ul{I}_T^0 : H^1(T) \to \ul{V}_T^0$ be the interpolator on $\ul{V}_T^0$ such that
\[
\ul{I}_T^0 v
\coloneq ( \lproj{F}{0} v )_{F \in \FT},
\]
where $\lproj{F}{0} : L^2(F) \to \Poly^0(F)$ is the $L^2$-orthogonal projector onto $\Poly^0(F)$ such that
\begin{equation}\label{eq:pi.P0}
  \int_F \lproj{F}{0} v = \int_F v.
\end{equation}
Let us check that
\begin{equation}\label{eq:pT.1:polynomial.consistency}
  p_T^1 (\ul{I}_T^0 v) = v \qquad \forall v \in \Poly^1(T).
\end{equation}
It holds, for all $w \in \Poly^1(T)$,
\[
\begin{aligned}
  \int_T \nabla p_T^1 (\ul{I}_T^0 v) \cdot \nabla w
  \overset{\eqref{eq:pT.1:gradient}}&= \sum_{F \in \FT} \int_F \lproj{F}{0} v (\nabla w \cdot n_{TF})
  \\
  \overset{\eqref{eq:pi.P0}}&= \sum_{F \in \FT} \int_F v (\nabla w \cdot n_{TF})
  = \int_T \nabla v \cdot \nabla w,
\end{aligned}
\]
where we have used the fact that $\nabla w \cdot n_{TF}$ is constant in the second equality, and an integration by parts to conclude.
This relation implies that
\begin{equation}\label{eq:nabla.pI=nabla}
  \nabla p_T^1 (\ul{I}_T^0 v) = \nabla v.
\end{equation}
On the other hand, \eqref{eq:pT.1:average} yields
\[
\int_T p_T^1 (\ul{I}_T^0 v)
= \sum_{F \in \FT} \frac{d_{TF}}{d} \int_F \lproj{F}{0} v
\overset{\eqref{eq:pi.P0}}=
\sum_{F \in \FT} \frac{d_{TF}}{d} \int_F v
= \int_T v,
\]
where the conclusion follows from the fact that $v$ is affine.
This relation shows that the average value of $p_T^1 (\ul{I}_T^0 v)$ inside $T$ is the same as that of $v$.
Combined with the equality of the gradients \eqref{eq:nabla.pI=nabla}, this yields \eqref{eq:pT.1:polynomial.consistency} after noticing that one can write, for all $x \in T$ and any $\phi\in\Poly^1(T)$, $\phi(x) = \frac{1}{|T|} \int_T \phi + \nabla \phi \cdot (x - \overline{x}_T)$.

\subsection{Extension to arbitrary-order}\label{sec:extension.higher.order}

The construction above provides an extension of the Crouzeix--Raviart space to general element geometries.
A legitimate question is whether this construction can be extended to an arbitrary order $k \ge 0$.

The key point is to devise an extension of the space in $\ul{V}_h^0$ with sufficient infomation to recover polynomial functions in $\Poly^{k+1}$ inside mesh elements.
Let $T \in \Th$.
We start by noticing that \eqref{eq:pT.1:gradient} mimicks the following integration by parts formula, valid for all $(v, w) \in H^1(T) \times C^\infty(\overline{T})$:
\[
\int_T \nabla v \cdot \nabla w
= - \int_T v \cdot \Delta w
+ \sum_{F \in \FT} \int_F v (\nabla w \cdot n_{TF}).
\]
When $v \in \Poly^{k+1}(T)$, to fully describe its gradient it is sufficient to write the above equation for $w \in \Poly^{k+1}(T)$.
In this case, $\Delta w \in \Poly^{k-1}(T)$ and $\nabla w \cdot n_{TF} \in \Poly^k(F)$ for all $F \in \FT$, and we can write
\begin{equation}\label{eq:reconstruct.nablav}
	\int_T \nabla v \cdot \nabla w
	= - \int_T \lproj{T}{k-1} v \cdot \Delta w
	+ \sum_{F \in \FT} \int_F \lproj{F}{k} v (\nabla w \cdot n_{TF}),
\end{equation}
where, for any $Y \in \Th \cup \Fh$ and any polynomial degree $\ell \ge 0$, the $L^2$-orthogonal projector $\lproj{Y}{\ell} : L^2(Y) \to \Poly^\ell(Y)$ is such that
\begin{equation}\label{eq:def.pil}
  \int_Y \lproj{Y}{\ell} v\, \phi = \int_Y v\, \phi
  \qquad \forall \phi \in \Poly^\ell(Y).
\end{equation}
Thus, when $v \in \Poly^{k+1}(T)$, \eqref{eq:reconstruct.nablav} shows that a full knowledge of $\nabla v$ is provided by $(\lproj{T}{k-1} v, (\lproj{F}{k} v)_{F \in \Fh})$.
On the other hand, we already know that the average value of $v$ in $T$ can be computed from the face values when $k = 0$ while, for $k \ge 1$, \eqref{eq:def.pil} with $Y=T$, $\ell=k-1\ge 0$ and $\phi \equiv 1 \in \Poly^0(T)\subset\Poly^{k-1}(T)$ shows that
\begin{equation}\label{eq:average.vT}
  \int_T v
  = \int_T \lproj{T}{k-1} v,
\end{equation}
i.e., the information concerning the average value of $v$ in $T$ is contained in $\lproj{T}{k-1} v$.

The above remark suggests the following arbitrary-order extension of the space $\ul{V}_h^0$, which is a variation of the Hybrid High-Order (HHO) space originally introduced in \cite{Di-Pietro.Ern.ea:14}:
For all integer $k \ge 0$,
\begin{equation}\label{eq:Vhk}
  \begin{aligned}
    \ul{V}_h^k \coloneq \big\{
    \ul{v}_h = \big( ( v_T )_{T \in \Th}, ( v_F )_{F \in \Fh} )\big)
    \,:\,
    &\text{$v_T \in \Poly^{k-1}(T)$ for all $T \in \Th$,}
    \\
    &\text{$v_F \in \Poly^k(F)$ for all $F \in \Fh$}
    \big\}.
  \end{aligned}
\end{equation}
The restriction of $\ul{V}_h^k$ to a mesh element $T \in \Th$ is obtained collecting the polynomial components on $T$ and its faces:
\[
\ul{V}_T^k \coloneq \big\{
\ul{v}_T = \big( v_T, ( v_F )_{F \in \FT} \big)
\,:\,
\text{$v_T \in \Poly^{k-1}(T)$
  and $v_F \in \Poly^k(F)$ for all $F \in \FT$}
\big\}.
\]
The definition of the potential reconstruction $p_T^{k+1} : \ul{V}_T^k \to \Poly^{k+1}(T)$ becomes:
For all $\ul{v}_T \in \ul{V}_T^k$,
\begin{equation}\label{eq:pT.k+1}
  \begin{gathered}
    \int_T \nabla p_T^{k+1} \ul{v}_T \cdot \nabla w
    = -\int_T v_T \Delta w
    + \sum_{F \in \FT} \int_F v_F (\nabla w \cdot n_{TF})
    \qquad \forall w \in \Poly^{k+1}(T),
    \\
    \int_T p_T^{k+1} \ul{v}_T =
    \begin{cases}
      \sum_{F \in \FT} \frac{d_{TF}}{d} \int_F v_F & \text{if $k = 0$},
      \\
      \int_T v_T & \text{if $k \ge 1$}.
    \end{cases}
  \end{gathered}
\end{equation}
For future use, we note the following equivalent reformulation of the first condition in~\eqref{eq:pT.k+1}, obtained integrating by parts the first term in the right-hand side:
\begin{multline}\label{eq:pT.k+1'}
  \int_T \nabla p_T^{k+1} \ul{v}_T \cdot \nabla w
  = -\int_T \nabla v_T \cdot \nabla w
  + \sum_{F \in \FT} \int_F (v_F - v_T) (\nabla w \cdot n_{TF})
  \\ \forall w \in \Poly^{k+1}(T).
\end{multline}

\subsection{Elliptic projector}

Define the local interpolator $\ul{I}_T^k : H^1(T) \to \ul{V}_T^k$ such that, for all $v \in H^1(T)$,
\[
\ul{I}_T^k \coloneq \big(
\lproj{T}{k-1} v, (\lproj{F}{k} v)_{F \in \FT}
\big)
\]
as well as the \emph{elliptic projector}
\begin{equation}\label{eq:elliptic.projector}
  \eproj{T}{k+1} \coloneq p_T^{k+1} \circ \ul{I}_T^k.
\end{equation}
Consider the special case where $v\in\Poly^{k+1}(T)$.
Writing the first condition in \eqref{eq:pT.k+1} for $\ul{v}_T = \ul{I}_T^k v$ and recalling \eqref{eq:reconstruct.nablav}, we infer that it holds, for all $w \in \Poly^{k+1}(T)$,
\begin{subequations}\label{eq:elliptic.projector:idempotency}
  \begin{equation}\label{eq:elliptic.projector:idempotency:gradient}
    \int_T \nabla \eproj{T}{k+1} v \cdot \nabla w
    = -\int_T \lproj{T}{k-1} v \Delta w
    + \sum_{F \in \FT} \int_F \lproj{F}{k} v (\nabla w \cdot n_{TF})
    = \int_T \nabla v \cdot \nabla w,
  \end{equation}
  so that $\nabla \eproj{T}{k+1} v = \nabla v$,since both $\eproj{T}{k+1}v$ and $v$ belong to $\Poly^{k+1}(T)$, and their gradient is therefore entirely determined by the left-hand and right-hand sides above.
  Moreover, by the second condition in \eqref{eq:pT.k+1},
  \begin{equation}\label{eq:elliptic.projector:idempotency:average}
    \int_T \eproj{T}{k+1} v
    = \begin{cases}
      \sum_{F \in \FT} \frac{d_{TF}}{d} \int_F v_F & \text{if $k = 0$},
      \\
      \int_T v_T & \text{if $k \ge 1$}
    \end{cases}
    \overset{\eqref{eq:cr:average.vT},\,\eqref{eq:average.vT}}= \int_T v,
  \end{equation}
\end{subequations}
showing that $\eproj{T}{k+1}v$ and $v$, which only differ by a constant since $\nabla \eproj{T}{k+1} v = \nabla v$, have the same average value over $T$, and thus that $\eproj{T}{k+1} v = v$.

Notice that, as a consequence of the previous discussion, $\eproj{T}{k+1}$ preserves functions in $\Poly^{k+1}(T)$, and is therefore a projector on $\Poly^{k+1}(T)$.
Under a mild mesh regularity assumption that essentially requires that each element $T \in \Th$ is star-shaped with respect to a ball of diameter uniformly comparable to $h_T$, it holds, for all $v \in H^{k+2}(T)$,
\begin{equation}\label{eq:elliptic.projector:approximation}
  \norm{\nabla(v - \eproj{T}{k+1} v)}{L^2(\partial T)^d}
  \lesssim h_T^{k+\frac12} \seminorm{v}{H^{k+2}(T)}.
\end{equation}

\begin{remark}[Elliptic projector]
  Let $v \in H^1(T)$ and set
  \[
  \overline{v}_T \coloneq \begin{cases}
    \frac{1}{|T|} \sum_{F \in \FT} \frac{d_{TF}}{d} \int_F v & \text{if $k = 0$},
    \\
    \frac{1}{|T|} \int_T v & \text{if $k \ge 1$}.
  \end{cases}
  \]
  The name elliptic projector stems from the fact that $\eproj{T}{k+1}$ satisfies
  \[
  \eproj{T}{k+1} v
  = \mathop{\operatorname{argmin}}\limits_{w \in \Poly^{k+1}(T),\, \int_T w = \int_T \overline{v}_T} \norm{\nabla(v - w)}{L^2(T)^d}^2,
  \]
  a problem strongly related to the Lagrangian energy for the elliptic equation \eqref{eq:poisson:strong}.
  The Euler equation for this problem reads
  \begin{equation}\label{eq:elliptic.projector:characterization}
    \begin{gathered}
      \int_T \nabla \eproj{T}{k+1} v \cdot \nabla w
      = \int_T \nabla v \cdot \nabla w \qquad \forall w \in \Poly^{k+1}(T),
      \\
      \int_T \eproj{T}{k+1} v = \int_T \overline{v}_T.
    \end{gathered}
  \end{equation}
\end{remark}


\section{A hybrid high-order scheme on general meshes}

\subsection{Discrete Poincar\'e inequality in hybrid spaces}

A key point to use Crouzeix--Raviart elements to discretize the Poisson problem was establishing a discrete version of the Poincar\'e inequality in Lemma~\ref{lem:poincare:cr}.
In this section, we prove a similar result for the HHO space \eqref{eq:Vhk}.
The differences with respect to the proof of Lemma~\ref{lem:poincare:cr} will provide an indication on how to design a numerical scheme.

As it will become clear in the proof of Lemma~\ref{lem:poincare:hho} (see, in particular, the first inequality in \eqref{eq:poincare:hho:key}), the natural discrete counterpart of the $H^1$-seminorm in the context of HHO spaces is the mapping $\norm{{\cdot}}{1,h} : \ul{V}_h^k \to \mathbb{R}^+$ such that, for all $\ul{v}_h \in \ul{V}_h^k$,
\begin{equation}\label{eq:norm.1h}
  \begin{gathered}
    \norm{\ul{v}_h}{1,h}^2
    \coloneq \sum_{T \in \Th} \norm{\ul{v}_T}{1,T}^2
    \\
    \text{
      with $\norm{\ul{v}_T}{1,T}^2
      \coloneq \norm{\nabla v_T}{L^2(T)^d}^2
      + h_T^{-1} \sum_{F \in \FT} \norm{v_F - v_T}{L^2(F)}^2$
      for all $T \in \Th$.
    }
  \end{gathered}
\end{equation}
Above and in what follows, to unify the treatment of the cases $k = 0$ (where there are no polynomials attached to mesh elements) and $k \ge 1$, it is understood that
\begin{equation}\label{eq:vT}
  v_T \coloneq \begin{cases}
    \frac{1}{|T|} \sum_{F \in \FT} \frac{d_{TF}}{d} \int_F v_F & \text{if $k = 0$},
    \\
    v_T & \text{if $k \ge 1$}.
  \end{cases}
\end{equation}
Given $\ul{v}_h \in \ul{V}_h^k$, we additionally let $v_h \in \Poly^{\max(0, k-1)}(\Th)$ be the (possibly discontinuous) piecewise polynomial function such that
\begin{equation}\label{eq:vh}
  (v_h)_{|T} \coloneq v_T \qquad \forall T \in \Th. 
\end{equation}

\begin{lemma}[Poincar\'e inequality in HHO spaces]\label{lem:poincare:hho}
  Let
  \[
  \ul{V}_{h,0}^k \coloneq \left\{
  \ul{v}_h \in \ul{V}_h^k \,:\,
  \text{$v_F = 0$ for all $F \in \Fhb$}
  \right\}.
  \]
  Then, there exists a real number $C_{\rm P} > 0$ independent of $h$, but possibly depending on the mesh regularity parameter and the polynomial degree $k$,
  such that, for any $\ul{v}_h \in \ul{V}_{h,0}^k$, it holds
  \begin{equation}\label{eq:poincare:hho}
    \norm{v_h}{L^2(\Omega)} \le C_{\rm P} \norm{\ul{v}_h}{1,h}.
  \end{equation}
\end{lemma}

\begin{proof}
  The key idea to prove the discrete Poincar\'e inequality for the Crouzeix--Raviart space was to express the projection on $\Poly^0(\Th)$ of $v_h$ as the divergence of a function in the global Raviart--Thomas space using the surjectivity of $\operatorname{div} : \RTN^1(\Th) \to \Poly^0(\Th)$.
  This idea, however, breaks down for HHO spaces because it would require to dispose of a Raviart--Thomas space on general polyhedral meshes.

  One way to circumvent this problem consists in noticing that $v_h \in L^2(\Omega)$ and in using instead the surjectivity of $\operatorname{div} : H^1(\Omega)^d \to L^2(\Omega)$ already invoked in the proof of Proposition~\ref{prop:surjectivity.div.RTN1.P0}.
  This gives the existence of $\tau \in H^1(\Omega)^d$ such that
  \begin{equation}\label{eq:surjectivity.div.H1.L2:bis}
    \text{%
      $\operatorname{div} \tau = v_h$ and
      $\norm{\tau}{H^1(\Omega)^d} \lesssim \norm{v_h}{L^2(\Omega)}$.
    }
  \end{equation}
  We can now continue in a similar way as we did in Lemma~\ref{lem:poincare:cr}:
  \[
  \begin{aligned}
    \norm{v_h}{L^2(\Omega)}^2
    \overset{\eqref{eq:surjectivity.div.H1.L2:bis}}&=
    \sum_{T \in \Th} \int_T v_T \operatorname{div} \tau
    \\
    &= \sum_{T \in \Th}\left[
      - \int_T \nabla v_T \cdot \tau
      + \sum_{F \in \FT} \int_F v_T (\tau \cdot n_{TF})
      \right]   
    \\
    \overset{\eqref{eq:magic.formula}}&= \sum_{T \in \Th}\left[
      - \int_T \nabla v_T \cdot \tau
      + \sum_{F \in \FT} \int_F (v_T - v_F) (\tau \cdot n_{TF})
      \right],
  \end{aligned}
  \]
  where we have used an integration by parts to pass to the second line.
  We can then apply Cauchy--Schwarz inequalities first to the integrals and then to the sums in the right-hand side to write
  \begin{align}
    \norm{v_h}{L^2(\Omega)}^2
    &\le \sum_{T \in \Th} \Bigg[
      \norm{\nabla v_T}{L^2(T)^d} \norm{\tau}{L^2(T)^d}
      \nonumber\\
      &\qquad\qquad
      + \sum_{F \in \FT} h_T^{-\frac12} \norm{v_F - v_T}{L^2(F)} ~
      h_T^{\frac12} \norm{\tau \cdot n_{TF}}{L^2(F)}
      \Bigg]
    \nonumber\\
    &\le \left[
      \sum_{T \in \Th}\left(
      \norm{\nabla v_T}{L^2(T)^d}^2
      + h_T^{-1} \sum_{F \in \FT} \norm{v_F - v_T}{L^2(F)}^2
      \right)
      \right]^{\frac12}
    \nonumber\\
    &\quad\times\left[
      \sum_{T \in \Th} \left(
      \norm{\tau}{L^2(T)^d}^2
      + h_T \norm{\tau}{L^2(\partial T)}^2        
      \right)
      \right]^{\frac12},
	  \label{eq:poincare:hho:key}
  \end{align}
  where, to pass to the last line, we have additionally used the H\"{o}lder inequality with exponents $(2,\infty,2)$ to write, for all $F \in \FT$,
  $\norm{\tau \cdot n_{TF}}{L^2(F)}
  \le \norm{n_{TF}}{L^\infty(F)^d} \norm{\tau}{L^2(F)^d}
  = \norm{\tau}{L^2(F)^d}$,
  the last equality following from the fact that $n_{TF}$ is a unit vector.
  The first factor in~\eqref{eq:poincare:hho:key} is precisely $\norm{\ul{v}_h}{1,h}$.
  The second factor can be bounded by $\norm{\tau}{H^1(\Omega)^d}$ using for the components of $\tau$ the local trace inequality, valid for all $T \in \Th$ and all $w \in H^1(T)$:
  \begin{equation}\label{eq:trace.cont}
    h_T^{\frac12}\norm{w}{L^2(\partial T)}
    \lesssim \norm{w}{L^2(T)}
    + h_T \norm{\nabla w}{L^2(T)^d}.
  \end{equation}
  We therefore obtain
  \[
  \norm{v_h}{L^2(\Omega)}^2
  \lesssim \norm{\ul{v}_h}{1,h} \norm{\tau}{H^1(\Omega)^d}
  \overset{\eqref{eq:surjectivity.div.H1.L2:bis}}\lesssim
  \norm{\ul{v}_h}{1,h} \norm{v_h}{L^2(\Omega)}.
  \]
  Simplifying, the conclusion follows.
\end{proof}

\begin{remark}[Norm $\norm{{\cdot}}{1,h}$]
  A simple consequence of Lemma~\ref{lem:poincare:hho} is that the mapping $\norm{{\cdot}}{1,h}$ is a norm on $\ul{V}_{h,0}^k$.
\end{remark}

\subsection{A coercive bilinear form}\label{sec:HHO:ah}

The HHO discretization of the Poisson problem hinges on the bilinear form $a_h : \ul{V}_h^k \times \ul{V}_h^k \to \mathbb{R}$ such that, for all $(\ul{w}_h, \ul{v}_h) \in \ul{V}_h^k \times \ul{V}_h^k$,
\begin{equation}\label{eq:ah.aT}
  \begin{gathered}
    a_h(\ul{w}_h, \ul{v}_h) \coloneq \sum_{T \in \Th} a_T(\ul{w}_T, \ul{v}_T)
    \\
    \text{
      with $a_T(\ul{w}_T, \ul{v}_T) \coloneq \int_T \nabla p_T^{k+1} \ul{w}_T \cdot \nabla p_T^{k+1} \ul{v}_T + s_T(\ul{w}_T, \ul{v}_T)$.
    }
  \end{gathered}
\end{equation}
The purpose of the symmetric positive semi-definite bilinear form $s_T : \ul{V}_T^k \times \ul{V}_T^k \to \mathbb{R}$ is to ensure coercivity and boundedness with respect to the $\norm{{\cdot}}{1,h}$-norm defined by \eqref{eq:norm.1h}, as made precise by the following assumption.

\begin{assumption}[Coercivity and boundedness]\label{ass:ST1}
  There is $\eta > 0$ independent of $h$ such that, for all $T \in \Th$,
  \begin{equation}\label{eq:ST1}\tag{ST1}
    \eta^{-1} \norm{\ul{v}_T}{1,T}^2
    \le a_T(\ul{v}_T, \ul{v}_T)
    \le \eta \norm{\ul{v}_T}{1,T}^2.
  \end{equation}
\end{assumption}

Summing \eqref{eq:ST1} over $T \in \Th$, we infer the coercivity and boundedness of the global bilinear form $a_h$:
\begin{equation}\label{eq:ah:coercivity.boundedness}
  \eta^{-1} \norm{\ul{v}_h}{1,h}^2
  \le a_h(\ul{v}_h, \ul{v}_h)
  \le \eta \norm{\ul{v}_h}{1,h}^2.
\end{equation}

\subsection{Discrete problem}

The discrete problem reads:
Find $\ul{u}_h \in \ul{V}_{h,0}^k$ such that
\begin{equation}\label{eq:poisson:discrete.basic}
  a_h(\ul{u}_h, \ul{v}_h) = \int_\Omega f v_h
  \qquad \forall \ul{v}_h \in \ul{V}_{h,0}^k,
\end{equation}
where $v_h \in \Poly^{\min(0, k-1)}(\Th)$ is the piecewise polynomial field defined by \eqref{eq:vh}.

\begin{remark}[Discretization of the right-hand side]
  Notice that the right-hand side is discretized using the function $v_h$ defined by~\eqref{eq:vh} and not the potential reconstruction.
  At this stage, this choice can be justified noticing that the norm of $v_h$ appears in the left-hand side of the discrete Poincaré inequality~\eqref{eq:poincare:hho}, and that this inequality is needed to prove stability reasoning as in Proposition~\ref{prop:galerkin:stability}; see Corollary~\ref{eq:poisson:weak}.
  This choice will also turn out to be crucial also to obtain consistency in Section~\ref{sec:basic.concepts:convergence.rate}.  
\end{remark}

\subsection{Stability}

The first step in the analysis of the HHO scheme is to obtain an a priori estimate on the discrete solution, from which uniqueness immediately follows.
Given a linear form $\ell_h : \ul{V}_{h,0}^k \to \mathbb{R}$, we define its dual norm setting
\begin{equation}\label{eq:dual.norm}
  \norm{\ell_h}{1,h,*} \coloneq
  \sup_{\ul{v}_h \in \ul{V}_{h,0}^k \setminus \{ \ul{0} \}} \frac{\ell_h(\ul{v}_h)}{\norm{\ul{v}_h}{1,h}}.
\end{equation}

\begin{lemma}[A priori estimate]\label{lem:a-priori}
  Denote by $\ell_h : \ul{V}_{h,0}^k \to \mathbb{R}$ a linear form.
  Let $\ul{w}_h \in \ul{V}_{h,0}^k$ be such that
  \begin{equation}\label{eq:poisson:auxiliary}
    a_h(\ul{w}_h, \ul{v}_h)
    = \ell_h(\ul{v}_h)
    \qquad\forall \ul{v}_h \in \ul{V}_{h,0}^k.
  \end{equation}
  Then, it holds
  \begin{equation}\label{eq:a-priori}
    \norm{\ul{w}_h}{1,h} \le \eta \norm{\ell_h}{1,h,*}.
  \end{equation}
\end{lemma}

\begin{proof}
  Taking $\ul{v}_h = \ul{w}_h$ in \eqref{eq:poisson:auxiliary} and using the coercivity of $a_h$ corresponding to the first inequality in \eqref{eq:ah:coercivity.boundedness}, we get
  \[
  \eta^{-1} \norm{\ul{w}_h}{1,h}^2
  \le a_h(\ul{w}_h, \ul{w}_h)
  = \ell_h(\ul{w}_h)
  \le \norm{\ell_h}{1,h,*} \norm{\ul{w}_h}{1,h},
  \]
  where the last passage follows from the definition \eqref{eq:dual.norm} of the dual norm.
  Simplifying by $\norm{\ul{w}_h}{1,h}$ and multiplying both sides by $\eta$ yields \eqref{eq:a-priori}.
\end{proof}

\begin{corollary}[Well-posedness of \eqref{eq:poisson:weak}]
  Problem \eqref{eq:poisson:weak} admits a unique solution that satisfies the a priori estimate
  \[
  \norm{\ul{u}_h}{1,h} \le \eta C_{\rm P} \norm{f}{L^2(\Omega)}.
  \]
\end{corollary}

\begin{proof}
  The existence and uniqueness of a solution follow observing that, once a basis for $\ul{V}_{h,0}^k$ is fixed, problem \eqref{eq:poisson:discrete.basic} is equivalent to an algebraic system with symmetric positive definite matrix.
  In order to prove the a priori bound on the solution, we apply Lemma~\ref{lem:a-priori} to the linear form $\ell_h : \ul{V}_{h,0}^k \ni \ul{v}_h \mapsto \int_\Omega f v_h \in \mathbb{R}$ and recall the definition \eqref{eq:dual.norm} of the dual norm to write
  \[
  \begin{aligned}
    \norm{\ul{u}_h}{1,h}
    &\le \eta \sup_{\ul{v}_h \in \ul{V}_{h,0}^k \setminus \{ \ul{0} \}} \frac{\int_\Omega f v_h}{\norm{\ul{v}_h}{1,h}}
    \\
    &\le \eta \sup_{\ul{v}_h \in \ul{V}_{h,0}^k \setminus \{ \ul{0} \}} \frac{\norm{f}{L^2(\Omega)} \norm{v_h}{L^2(\Omega)}}{\norm{\ul{v}_h}{1,h}}
    \\
    \overset{\eqref{eq:poincare:hho}}&\le \eta C_{\rm P}
    \sup_{\ul{v}_h \in \ul{V}_{h,0}^k \setminus \{ \ul{0} \}} \frac{\norm{f}{L^2(\Omega)} \cancel{\norm{\ul{v}_h}{1,h}}}{\cancel{\norm{\ul{v}_h}{1,h}}}
    = \eta C_{\rm P} \norm{f}{L^2(\Omega)}. \qedhere
  \end{aligned}
  \]
\end{proof}

\subsection{Energy error estimate}

To study the convergence of the scheme \eqref{eq:poisson:discrete.basic}, we must identify a suitable measure of the error.
This question is not entirely trivial, as the exact solution $u$ and its numerical approximation $\ul{u}_h$ are not directly comparable: we cannot identify a suitable function space to which they both belong.
Two choices are possible: either reconstructing from $\ul{u}_h$ a function in, e.g., $H_0^1(\Omega)$ and compare it with $u$, or reducing $u$ to an element of $\ul{V}_{h,0}^k$ and compare the latter with $\ul{u}_h$.
It was argued in \cite{Di-Pietro.Droniou:18} that the latter choice leads to a more natural analysis, where
each of the classical properties (stability, consistency, boundedness) leads to a well-identified major consequence and
fewer terms have to be estimated to prove consistency.
We therefore define the error
\begin{equation}\label{eq:eh}
  \ul{e}_h \coloneq \ul{u}_h - \ul{I}_h^k u \in \ul{V}_{h,0}^k
\end{equation}
where $\ul{I}_h^k : H^1(\Omega) \to \ul{V}_h^k$ is the global interpolator such that, for all $v \in H^1(\Omega)$,
\[
\ul{I}_h^k v \coloneq \big(
( \lproj{T}{k-1} v)_{T \in \Th},
( \lproj{F}{k} v)_{F \in \Fh}
\big).
\]
Notice that, when $u\in H^1_0(\Omega)$, $\ul{I}_h^k u\in \ul{V}_{h,0}^k$.

The key idea to bound $\norm{\ul{e}_h}{1,h}$ in terms of a consistency error is to use the a priori estimate \eqref{eq:a-priori} for the following problem:
For all $\ul{v}_h \in \ul{V}_{h,0}^k$,
\begin{multline}\label{eq:err.eq}
  a_h(\ul{e}_h, \ul{v}_h)
  \overset{\eqref{eq:eh}}=
  a_h(\ul{u}_h, \ul{e}_h) - a_h(\ul{I}_h^k u, \ul{v}_h)
  \\
  \overset{\eqref{eq:poisson:discrete.basic}}=
  \int_\Omega f v_h - a_h(\ul{I}_h^k u, \ul{v}_h)
  \eqcolon \mathcal{E}_h(u; \ul{v}_h),
\end{multline}
which is precisely the equation for the error.
We thus have the basic error estimate
\begin{equation}\label{eq:basic.error.estimate}
  \norm{\ul{e}_h}{1,h} \le \eta \norm{\mathcal{E}_h(u; \cdot)}{1,h,*}.
\end{equation}
This estimate shows that the method is convergent if the following consistency property holds:
\begin{equation}\label{eq:conv.consistency}
	\lim_{h \to 0^+} \norm{\mathcal{E}_h(u; \cdot)}{1,h,*} = 0.
\end{equation}
This is true if, in particular, $\norm{\mathcal{E}_h(u; \cdot)}{1,h,*} \lesssim h^{\alpha}$ with $\alpha > 0$ and hidden constant independent of $h$ (but possibly depending on $\Omega$, $d$, the mesh regularity parameter, and the polynomial degree $k$).

\begin{remark}[Optimality of the basic error estimate]
  From the error equation $\mathcal{E}_h(u; \ul{v}_h)=a_h(\ul{e}_h, \ul{v}_h)$, a Cauchy--Schwarz inequality, and the global boundedness \eqref{eq:ah:coercivity.boundedness} of $a_h$, we easily obtain the following converse inequality of \eqref{eq:basic.error.estimate}:
  \begin{equation}\label{eq:optimality}
    \norm{\mathcal{E}_h(u; \cdot)}{1,h,*}\le \eta \norm{\ul{e}_h}{1,h},
  \end{equation}
  showing that the convergence \eqref{eq:conv.consistency} towards $0$ of $\mathcal{E}_h(u; \cdot)$ is actually also a necessary condition to the convergence of the error $\ul{e}_h$ towards $0$.
  Additionally, since the multiplicative constant $\eta$ in \eqref{eq:basic.error.estimate} and \eqref{eq:optimality} does not depend on $h$, if we prove a convergence rate in $h$ for $\norm{\mathcal{E}_h(u; \cdot)}{1,h,*}$ (which will make the object of the following section), then this convergence rate will be the optimal one for $\norm{\ul{e}_h}{1,h}$.
\end{remark}

\subsection{Convergence rate of the consistency error for smooth solutions}\label{sec:basic.concepts:convergence.rate}

A key step to estimate the consistency error consists in reformulating it so that differences of the form $u - \eproj{T}{k+1} u$ appear, and then use the approximation properties of the elliptic projector to identify the order of convergence in $h$ for smooth solution.
We will assume throughout the rest of this section that $u$ is sufficiently regular inside each element to fully exploit the approximation properties of the HHO space of degree $k$, i.e.,
\[
u \in H^{k+2}(\Th) \coloneq \left\{
v \in L^2(\Omega) \,:\, \text{$v_{|T} \in H^{k+2}(T)$ for all $T \in \Th$}
\right\},
\]
so that
\[
\seminorm{u}{H^{k+2}(\Th)} \coloneq \left(
\sum_{T \in \Th} \seminorm{u}{H^{k+2}(T)}^2
\right)^{\frac12} < +\infty.
\]

Let $\ul{v}_h \in \ul{V}_{h,0}^k$ and consider each of the contributions to $\mathcal{E}_h(u; \ul{v}_h)$ (see~\eqref{eq:err.eq}) separately.
For the first one, we recall that $f = - \Delta u$ almost everywhere in $\Omega$ to write
\[
\begin{aligned}
  \int_\Omega f v_h
  &= - \int_\Omega \Delta u \, v_h
  = - \sum_{T \in \Th} \int_T \Delta u \, v_T
  \\
  &= \sum_{T \in \Th} \left(
  \int_T \nabla u \cdot \nabla v_T
  - \sum_{F \in \FT} \int_F \nabla u \cdot n_{TF} v_T
  \right),
\end{aligned}
\]
where we have integrated by parts element by element in the last passage.
We next notice that $\nabla u \in H(\operatorname{div}; \Omega)$ (since $\nabla u \in L^2(\Omega)^d$ as a consequence of $u \in H_0^1(\Omega)$ and $\operatorname{div} \nabla u = \Delta u \overset{\eqref{eq:poisson:strong}}= -f \in L^2(\Omega)$) to justify the use of Proposition~\ref{prop:magic.formula} and insert $v_F$ into the boundary term, thus arriving at
\begin{equation}\label{eq:Eh:T1}
  \int_\Omega f v_h
  = \sum_{T \in \Th} \left(
  \int_T \nabla u \cdot \nabla v_T
  + \sum_{F \in \FT} \int_F \nabla u \cdot n_{TF} (v_F - v_T)
  \right).
\end{equation}

Let us consider the second contribution in $\mathcal{E}_h(u; \ul{v}_h)$.
Using the definition \eqref{eq:ah.aT} of $a_h$ and $a_T$, and recalling the definition \eqref{eq:elliptic.projector} of the elliptic projector, we have
\begin{equation}\label{eq:Eh:T2}
\begin{aligned}
  a_h(\ul{I}_h^k u, \ul{v}_h)
  &= \sum_{T \in \Th} \left(
  \int_T \nabla \eproj{T}{k+1} u \cdot \nabla p_T^{k+1} \ul{v}_T
  + s_T(\ul{I}_T^k u, \ul{v}_T)
  \right)
  \\
  \overset{\eqref{eq:pT.k+1'}}&=
  \sum_{T \in \Th} \left(
  \int_T \nabla \eproj{T}{k+1} u \cdot \nabla v_T
  + \sum_{F \in \FT} \int_F \nabla \eproj{T}{k+1} u \cdot n_{TF} (v_F - v_T)
  \right)    
  \\
  &\quad
  + \sum_{T \in \Th} s_T(\ul{I}_T^k u, \ul{v}_T),  
\end{aligned}
\end{equation}
Subtracting \eqref{eq:Eh:T2} from \eqref{eq:Eh:T1}, we arrive at the following reformulation of the consistency error:
\begin{equation}\label{eq:Eh:decomposition}
  \begin{aligned}
    \mathcal{E}_h(u; \ul{v}_h)
    &=
    \underbrace{%
      \sum_{T \in \Th} \int_T \nabla (u - \eproj{T}{k+1} u) \cdot \nabla v_T
    }_{\term_1}
    \\
    &\quad
    + \underbrace{%
      \sum_{T \in \Th} \sum_{F \in \FT} \int_F \nabla (u - \eproj{T}{k+1} u) \cdot n_{TF} (v_F - v_T)
    }_{\term_2}
    - \underbrace{%
      \sum_{T \in \Th} s_T(\ul{I}_T^k u, \ul{v}_T).
    }_{\term_3}
  \end{aligned}
\end{equation}
We proceed to estimate the terms in the right-hand side.
Recalling the characterization \eqref{eq:elliptic.projector:characterization} of the elliptic projector and noticing that $v_T \in \Poly^{k-1}(T) \subset \Poly^{k+1}(T)$, it is readily inferred that
\begin{equation}\label{eq:Eh:estimate.T1}
  \term_1 = 0.
\end{equation}

For the second term, we use the H\"older inequality with exponents $(2,\infty,2)$ on the integral to write
\begin{equation}\label{eq:Eh:estimate.T2}
  \begin{aligned}
    \term_2
    &\le
    \sum_{T \in \Th} \sum_{F \in \FT}
    h_T^{\frac12}\norm{\nabla (u - \eproj{T}{k+1} u)}{L^2(F)^d} \,
    \norm{n_{TF}}{L^\infty(F)^d} \,
    h_T^{-\frac12}\norm{v_F - v_T}{L^2(F)}
    \\
    &\le
    \left(
    \sum_{T \in \Th} h_T \norm{\nabla(u - \eproj{T}{k+1} u)}{L^2(\partial T)^d}^2
    \right)^{\frac12} \left(
    \sum_{T \in \Th} \sum_{F \in \FT} h_T^{-1} \norm{v_F - v_T}{L^2(F)}^2
    \right)^{\frac12}
    \\
    \overset{\eqref{eq:elliptic.projector:approximation},\,\eqref{eq:norm.1h}}&\lesssim
    h^{k+1} \seminorm{u}{H^{k+2}(\Th)} \norm{\ul{v}_h}{1,h},
  \end{aligned}
\end{equation}
where, to pass to the second line, we have used a Cauchy--Schwarz inequality on the sums, along with the fact that $\norm{n_{TF}}{L^\infty(F)^d} = 1$
while, in the conclusion, we have additionally used the definition \eqref{eq:norm.1h} of $\norm{\ul{v}_h}{1,h}$ together with the fact that $h_T \le h$ for all $T \in \Th$.

It only remains to estimate the term involving the stabilization.
To this end, we will need an additional assumption on the stabilization bilinear form.

\begin{assumption}[Polynomial consistency]\label{ass:ST2}
  For all $T \in \Th$ and all $w \in \Poly^{k+1}(T)$,
  \begin{equation}\label{eq:ST2}\tag{ST2}
    s_T(\ul{I}_T^k w, \ul{v}_T) = 0
    \qquad \forall \ul{v}_T \in \ul{V}_T^k.
  \end{equation}
\end{assumption}

In the next section we will show that \eqref{eq:ST1} and \eqref{eq:ST2} imply
\begin{equation}\label{eq:Eh:estimate.T3}
  \term_3 \lesssim h^{k+1} \seminorm{u}{H^{k+2}(\Th)} \norm{\ul{v}_h}{1,h}.
\end{equation}

Using the estimates \eqref{eq:Eh:estimate.T1}, \eqref{eq:Eh:estimate.T2}, and \eqref{eq:Eh:estimate.T3} to bound the right-hand side of \eqref{eq:Eh:decomposition}, we conclude that $\mathcal{E}_h(u; \ul{v}_h) \lesssim h^{k+1} \seminorm{u}{H^{k+2}(\Th)} \norm{\ul{v}_h}{1,h}$. Dividing by $\norm{\ul{v}_h}{1,h}$ and passing to the supremum over $\ul{v}_h \in \ul{V}_{h,0}^k \setminus \{ \ul{0} \}$, this gives
\begin{equation}\label{eq:Eh:estimate}
  \norm{\mathcal{E}_h(u; \cdot)}{1,h,*} \lesssim h^{k+1} \seminorm{u}{H^{k+2}(\Th)}.
\end{equation}

We conclude with a theorem that summarises the convergence rate of the error estimate that we have obtained.

\begin{theorem}[Convergence rate]\label{thm:convergence.rate}
  Let $u \in H_0^1(\Omega)$ and $\ul{u}_h \in \ul{V}_{h,0}^k$ solve \eqref{eq:poisson:weak} and \eqref{eq:poisson:discrete.basic}, respectively.
  Further assume that $u \in H^{k+2}(\Th)$ and that, for all $T \in \Th$, $s_T$ satisfies Assumptions~\ref{ass:ST1} and \ref{ass:ST2}.
  Then,
  \[
  \norm{\ul{I}_h^ku - \ul{u}_h}{1,h} \lesssim h^{k+1} \seminorm{u}{H^{k+2}(\Th)}.
  \]
\end{theorem}

\begin{proof}
  It suffices to plug \eqref{eq:Eh:estimate} into \eqref{eq:basic.error.estimate}.
\end{proof}

\subsection{Consistency of the stabilization for smooth functions}

The goal of this section is to prove the following lemma, from which \eqref{eq:Eh:estimate.T3} immediately follows.
The proof hinges on two results described hereafter.
The first result is the $H^s$-boundedness of the $L^2$-orthogonal projector:
For all $T \in \Th$ and all integers $\ell \ge 0$ and $s \ge 0$, it holds
\begin{equation}\label{eq:l2proj:boundedness}
  \seminorm{\lproj{T}{\ell} v}{H^s(T)} \le \seminorm{v}{H^s(T)}
  \qquad \forall v \in H^s(T),
\end{equation}
with hidden constant independent of the meshsize but possibly depending on $d$, $\ell$, $s$, and the mesh regularity parameter.
Using~\eqref{eq:l2proj:boundedness} in conjunction with the ideas of~\cite{Dupont.Scott:80}, the following approximation properties for $\lproj{T}{\ell}$ on elements of general shape were proved in~\cite{Di-Pietro.Droniou:17} (see also~\cite[Chapter~1]{Di-Pietro.Droniou:20}):
For all $T \in \Th$,
all $v \in H^s(T)$ with $s \in \{ 0, \ldots, \ell + 1 \}$,
and all $m \in \{ 0, \ldots, s \}$,
\begin{subequations}\label{eq:polynomial.approximation}
  \begin{equation}\label{eq:polynomial.approximation:element}
    \norm{\partial^{\alpha} (v - \lproj{T}{\ell} v)}{L^2(T)}
    \lesssim h_T^{s-m} \seminorm{v}{H^s(T)},
  \end{equation}
  and, if $s \ge 1$ and $m \le s - 1$,
  \begin{equation}\label{eq:polynomial.approximation:trace}
    \norm{\partial^{\alpha} (v - \lproj{T}{\ell} v)}{L^2(\partial T)}
    \lesssim h_T^{s-m-\frac12} \seminorm{v}{H^s(T)},
  \end{equation}
\end{subequations}
where $\alpha \in \mathbb{N}^d$ is a vector of positive integers such that $\sum_{i=1}^d \alpha_i \le m$ and $\partial^{\alpha} \coloneq \partial_{\alpha_1} \cdots \partial_{\alpha_d}$.
The hidden constants in \eqref{eq:polynomial.approximation} depend only on $d$, $l$, $s$, and the mesh regularity parameter, but are independent of the mesh size.

\begin{lemma}[Consistency of the stabilization for smooth functions]\label{lem:sT:consistency}
  Let $T \in \Th$ and let $s_T : \ul{V}_T^k \times \ul{V}_T^k \to \mathbb{R}$ be a symmetric positive semi-definite bilinear form satisfying Assumptions~\ref{ass:ST1} and \ref{ass:ST2}.
  Then, for all $v \in H^{k+2}(T)$, it holds
  \begin{equation}\label{eq:sT:consistency}
    s_T(\ul{I}_T^k v, \ul{I}_T^k v)^{\frac12}
    \lesssim h_T^{k+1} \seminorm{v}{H^{k+2}(T)},
  \end{equation}
  with hidden constant independent of $h$ and $T$, but possibly depending on $d$, on the mesh regularity parameter, and on the polynomial degree $k$.
\end{lemma}

\begin{proof}
  By the bilinearity and symmetry of $s_T$ followed by \eqref{eq:ST2}, we have
  \[
  \begin{aligned}
    &s_T(\ul{I}_T^k (v - \lproj{T}{k+1} v), \ul{I}_T^k (v - \lproj{T}{k+1} v))
    \\
    &\quad=
    s_T(\ul{I}_T^k v, \ul{I}_T^k v)
    -2s_T(\ul{I}_T^k v, \ul{I}_T^k  \lproj{T}{k+1} v)
    +s_T(\ul{I}_T^k  \lproj{T}{k+1} v, \ul{I}_T^k \lproj{T}{k+1} v)
    \\
    &\quad=
    s_T(\ul{I}_T^k v, \ul{I}_T^k v).
  \end{aligned}
  \]
  Hence,
  \begin{equation}\label{eq:sT:consistency:basic}
    \begin{aligned}
      s_T(\ul{I}_T^k v, \ul{I}_T^k v)
      &= s_T(\ul{I}_T^k (v - \lproj{T}{k+1} v), \ul{I}_T^k (v - \lproj{T}{k+1} v))
      \\
      \overset{\eqref{eq:ST1}}&\lesssim
      \norm{\ul{I}_T^k (v - \lproj{T}{k+1} v)}{1,T}^2.
    \end{aligned}
  \end{equation}
  By definition \eqref{eq:norm.1h} of the $\norm{{\cdot}}{1,T}$-norm, we have
  \[
  \begin{aligned}
    &\norm{\ul{I}_T^k (v - \lproj{T}{k+1} v)}{1,T}^2
    \\
    &\quad=
    \norm{\nabla \lproj{T}{k-1} (v - \lproj{T}{k+1} v)}{L^2(T)^d}^2
    + h_T^{-1} \sum_{F \in \FT} \norm{\lproj{F}{k} (v - \lproj{T}{k+1} v) - \lproj{T}{k-1} (v - \lproj{T}{k+1} v)}{L^2(F)}^2
    \\
    &\quad\lesssim
    \norm{\nabla \lproj{T}{k-1} (v - \lproj{T}{k+1} v)}{L^2(T)^d}^2
    + h_T^{-1} \sum_{F \in \FT} \norm{\lproj{F}{k} (v - \lproj{T}{k+1} v)}{L^2(F)}^2
    \\
    &\qquad
    + h_T^{-1} \norm{\lproj{T}{k-1} (v - \lproj{T}{k+1} v)}{L^2(\partial T)}^2
    \\
    &\quad\eqcolon \term_1 + \term_2 + \term_3,
  \end{aligned}
  \]
  where we have used a triangle inequality together with $(a + b)^2 \le 2(a^2 + b^2)$ on the boundary term in the second passage.
  Using the boundedness property \eqref{eq:l2proj:boundedness} of the $L^2$-orthogonal projector with $(\ell, s) = ( k-1, 1)$, we infer that
  \[
  \term_1 \lesssim \norm{\nabla (v - \lproj{T}{k+1} v)}{L^2(T)^d}^2.
  \]
  Noticing that $\lproj{F}{k}$ is a contraction mapping for the $L^2(F)$-norm, we can bound
  \[
  \sum_{F \in \FT} \norm{\lproj{F}{k} (v - \lproj{T}{k+1} v)}{L^2(F)}^2
  \lesssim h_T^{-1} \norm{v - \lproj{T}{k+1} v}{L^2(\partial T)}^2
  \lesssim h_T^{2(k+1)} \seminorm{v}{H^{k+2}(T)}^2,
  \]
  where the conclusion follows from \eqref{eq:polynomial.approximation:trace} with $(\ell, s, m) = (k + 1, k + 2, 0)$.
  Finally, applying the trace inequality~\eqref{eq:trace.cont} to $w = v - \lproj{T}{k+1} v$, we infer for the last term
  \[
  \begin{aligned}
  \term_3
  &\lesssim
  h_T^{-2} \norm{\lproj{T}{k-1} (v - \lproj{T}{k+1} v)}{L^2(T)}^2
  + \norm{\nabla \lproj{T}{k-1} (v - \lproj{T}{k+1} v)}{L^2(T)^d}^2
  \\
  \overset{\eqref{eq:l2proj:boundedness}}&\lesssim
  h_T^{-2} \norm{v - \lproj{T}{k+1} v}{L^2(T)}^2
  + \norm{\nabla (v - \lproj{T}{k+1} v)}{L^2(T)^d}^2
  \\
  &\lesssim h_T^{2(k+1)} \seminorm{v}{H^{k+2}(T)}^2,
  \end{aligned}
  \]
  where the conclusion follows from \eqref{eq:polynomial.approximation} with $\ell = k+1$, $s = k+2$, and $m$ respectively equal to $0$ and $1$.
  Gathering the above bounds for $\term_1$, $\term_2$, and $\term_3$, we conclude that
  \[
  \norm{\ul{I}_T^k (v - \lproj{T}{k+1} v)}{1,T}^2
  \lesssim h^{2(k+1)} \seminorm{v}{H^{k+2}(T)}^2.
  \]
  Plugging this bound in \eqref{eq:sT:consistency:basic} yields \eqref{eq:sT:consistency}.
\end{proof}

\subsection{An example of stabilization bilinear form}

To close this chapter, we provide an example of stabilization bilinear form that satisfies properties \eqref{eq:ST1}--\eqref{eq:ST2}.
Given $T \in \Th$, we consider the following expression:
For all $(\ul{w}_T, \ul{v}_T) \in \ul{V}_T^k \times \ul{V}_T^k$,
\[
s_T(\ul{w}_T, \ul{v}_T)
\coloneq
h_T^{-2} \int_T \delta_T^{k-1} \ul{w}_T \, \delta_T^{k-1} \ul{v}_T
+ h_T^{-1} \sum_{F \in \FT} \int_F \delta_{TF}^k \ul{w}_T \, \delta_{TF}^k \ul{v}_T,
\]
where
\[
\text{
  $\delta_T^{k-1} \ul{v}_T \coloneq v_T - \lproj{T}{k-1} (p_T^{k+1} \ul{v}_T)$
  and $\delta_{TF}^k \ul{v}_T \coloneq v_F - \lproj{F}{k} (p_T^{k+1} \ul{v}_T)$
  for all $F \in \FT$.
}
\]
In other words, $s_T$ penalizes in a least square sense the differences between the vector of polynomials $\ul{v}_T$ and $\ul{I}_T^k p_T^{k+1} \ul{v}_T$, the interpolate of the potential reconstruction obtained from the latter.

\begin{remark}[Link with the Crouzeix--Raviart scheme]
  If $\Th$ is a conforming simplicial mesh and $k = 0$, it holds $\delta_T^{k-1} \ul{v}_T \equiv 0$ since $k-1 = -1$ and $\Poly^{-1}(T) = \{ 0 \}$, and, for all $F \in \FT$, $\delta_{TF}^0 \ul{v}_T \equiv 0$ since $p_T^1 \ul{v}_T$ is the Crouzeix--Raviart function with degrees of freedom $v_F$, so that $\lproj{F}{0}(p_T^1 \ul{v}_T) = v_F$ for all $F \in \FT$.
  As a consequence, at the algebraic level the left-hand sides of problems~\eqref{eq:poisson:discrete.basic} and~\eqref{eq:cr:discrete} coincide.
  This is not the case, however, for the right-hand sides: the HHO scheme actually corresponds to replacing $\int_\Omega f v_h$ with $\sum_{T \in \Th} \int_T f\, \lproj{T}{0} v_h = \sum_{T \in \Th} \int_T \lproj{T}{0} f\, v_h$ in~\eqref{eq:cr:discrete}, a variation leading to the same orders of convergence as the original Crouzeix--Raviart scheme.
\end{remark}


\section*{Acknowledgements}

Funded by the European Union (ERC Synergy, NEMESIS, project number 101115663).
Views and opinions expressed are however those of the authors only and do not necessarily reflect those of the European Union or the European Research Council Executive Agency. Neither the European Union nor the granting authority can be held responsible for them.


\printbibliography

\end{document}